\documentclass[a4paper]{article} %
\usepackage{graphicx,amssymb} %
\date{} %

\usepackage{amsmath,amsthm}
\usepackage{amssymb}
\usepackage[T1]{fontenc}
\usepackage{amsmath}
\usepackage{amsfonts}
\usepackage{amssymb}
\usepackage{amsthm}
\usepackage{graphicx}
\usepackage{makeidx}
\usepackage{color}

\newcommand{\R}{\mathbb{R}}
\newcommand{\N}{\mathbb{N}}

\newtheorem{thm}{Theorem}[section]
\newtheorem{cor}[thm]{Corollary}
\newtheorem{lem}[thm]{Lemma}
\newtheorem{prop}[thm]{Proposition}

\newtheorem{defn}[thm]{Definition}

\theoremstyle{definition}
\newtheorem{defin}[thm]{Definition}
\newtheorem{rem}[thm]{Remark}
\newtheorem{exa}[thm]{Example}
\newtheorem{pre}[thm]{Proof}

\numberwithin{equation}{section}

\usepackage{enumerate}
\begin{document}


\baselineskip=17pt


\title{Some non noetherian $C^\infty$\\ quasianalytic local rings}

\author{Abdelhafed Elkhadiri\\
University Ibn Tofail,
Faculty of Sciences\\
Kenitra, Morocco\\
E-mail: elkhadiri.abdelhafed@uit.ac.ma
}

\date{}

\maketitle


\renewcommand{\thefootnote}{}



\renewcommand{\thefootnote}{\arabic{footnote}}


\begin{abstract}
We give an example of a non-noetherian quasi-analytic ring constructed using a quasi-analytic Denjoy-Carleman class.
\\ If we denote   by  $ \mathcal{D}_n$  the ring of those $ C^\infty$ quasianalytic function germs at $0\in \R^n$ which are definable in a polynomially bounded o-minimal structure.
We show that  the  system \\$\{ \mathcal{D}_n\,/\, n\in\N^*\}$ is not noetherian, i.e. there exists $m\in\N $, $ m > 1$, such that the ring $\mathcal{D}_m$ is not noetherian.
\end{abstract}

\section{Introduction}
Recall that a ring $\mathcal{C}_n$ of smooth germs at the origin of $\R^n$ is called quasianalytic if the only element of $\mathcal{C}_n$ that admits a zero Taylor expansion is the zero germ. Let, for each $n\in\N$, a quasianalytic ring $\mathcal{C}_n$, we say that 
$\mathcal{C}=\{ \mathcal{C}_n\,/\,n\in\N\}$ is a quasianalytic system, if $\mathcal{C}$ is closed under composition, partial differentiation and implicit function. \\ Hadamard proposed the following problem to to characterize quasianalytic rings:\\
Give necessary and sufficient conditions bearing on the growth of partial derivatives of $C^\infty$ germs at the origin, in order that the ring $\mathcal{C}_n$ is quasianalytic. \\
Denjoy \cite{Denjoy} gave conditions sufficient for quasianalyticity. But the problem was solved completely by Carleman \cite{Carleman}. This lead to the notion of quasianalytic Denjoy-Carleman classes of functions, see section 3 below. The algebraic properties of such rings, namely their stability under several classical operations, such as composition, differentiation, implicit function, is well understood see \cite{Childress}, \cite{Komatsu}.\\
However, two classical properties, namely Weierstrass division and Weierstrass preparation, are not valid in the quasianalytic sitting. It has been proved by Childress \cite{Childress} that quasianalytic Denjoy-Carleman classes might not satisfy Weierstrass division. On the other hand, Weierstrass preparation is also not valid
in the quasianalytic sitting, see \cite{Parusinski}, \cite{Broglia}. Because of the lack of Weierstrass division in quasianalytic Denjoy-Carleman classes many problems remain open in dimension $> 1$, for instance these classes are not known to be noetherian and to satisfy any kind of M. Artin approximation theorem. Moreover, we do not know any example of a non noetherian quasianalytic ring of germs, at the origin in $\R^n$, $ n >1$,  or a non noetherian quasianalytic Denjoy-Carleman class.\\
In this paper we give two examples of non noetherian quasianalytic rings.\\
- The first example is obtained as follow: we give a quasianalytic Denjoy-Carleman class such that all its shifted Denjoy-Carleman classes are quasianalytic, see section 5. We consider the union of all those classes. Using a result shown in \cite{Elkhadiri}, we prove that this ring is not noetherian.\\
- The second example arise from model theory.\\
We say that a quasianalytic system $ \mathcal{C} = \{ \mathcal{C}_n\,/\, n\in\N\}$ is noetherian if all the rings $\mathcal{C}_n$ are noetherian. If $\mathcal{R}$ is a polynomially bounded o-minimal structure, we denote by $ \mathcal{D}_n\subset \mathcal{E}_n$ the ring of germs, at the origin in $\R^n$, of smooth definable functions in the structure $\mathcal{R}$. By \cite{Miller}, the ring $ \mathcal{D}_n$ is quasianalytic. If the system $ \mathcal{D} = \{ \mathcal{D}_n\,/\, n\in\N\}$  contains strictly the  analytic system, we prove that the system $\{ \mathcal{D}_n\,/\, n\in\N\}$ is not noetherian.\\
Finally the definition that we give for  functions in a Quasianalytic Denjoy-Carleman class in section 3, is slightly different from the one usually given in the literature. It seems to us that it is simpler. But the reader may note that these two definitions are the same.
\section{Notations and Definition} 
Let $I$ denote the interval $ [0,1]\subset \R$ and for $ \epsilon > 0$, $ I_\epsilon = [0, \epsilon]$.
We denote by $\mathcal{E}(I^n)$ the ring of $C^\infty$ functions on 
$I^n = \underbrace{I\times  \ldots \times I}_{n \mbox{ times}}$. 
  $\mathcal{E}_n$  denotes the ring of smooth germs at the origin in $\R^n$ and  $ \R[[X_1,\ldots, X_n]]$ is the ring of formal power series with coefficients in $\R$.
For every $ f\in\mathcal{E}_n$, we denote by $\hat{f}\in \R[[X_1,\ldots, X_n]]$  its (infinite) Taylor expansion at the origin. The mapping $$ \mathcal{E}_n\ni f \mapsto \hat{f}\in \R[[X_1,\ldots, X_n]]$$ is called the Borel mapping.\\
Let $ n\in\N$, $ \alpha=(\alpha_1,\ldots,\alpha_n)\in\N^n$ and $ x= (x_1,\ldots,x_n)$, we use the standard notations
\[ |\alpha| = \alpha_1 + \ldots + \alpha_n,\,\,\,\,\alpha! = \alpha_1!\ldots \alpha_n!, \,\,\,\, D^\alpha = \frac{\partial^{|\alpha|}}
{\partial x_1^{\alpha_1}\ldots \partial x_n^{\alpha_n}}.\]
We say that a real function, $m$, of one real variable is defined [resp.$ C^\infty$] for $ t \gg 0 $, if there exists $ b > 0$ such that the function $m$ is defined [resp.$ C^\infty$] on the interval $ [ b,  \infty[ $. \\
In all the following, $m$ will be a $ C^\infty$ for $ t \gg 0$ such that:
\begin{enumerate}
\item for $ t \gg 0$, $m(t) >0$, $m'(t) >0$, $m''(t) >0$,
\item $\lim\limits_{t\to \infty} m'(t) = \infty$,
\item there exists $ \delta > 0$ such that $ m''(t) \leq \delta $ for $ t \gg 0$.
\end{enumerate}
We put $$ M(t)= e^{m(t)}.$$
If a function $ m:[b,\infty[\to \R$ is such that $m(b)=0$, we still denote by $m$ the extension of $m$ to $ [0,\infty[$ obtained by sitting $ m(t)=0$ if $ 0\leq t\leq b$. We note that this extension remains convex.
\section{Functions of the class $M$.}
\begin{defn}
A function $ f\in\mathcal{E}(I^n)$ is said to be in the class $M$, if there exist $ C >0$, $\rho >0$ such that
\[ \forall x\in I^n,\,\,\,\,|D^\alpha f(x)| \leq C \rho^{|\alpha|}M(|\alpha|),\,\,\,\,\forall \alpha\in\N^n,\,\,|\alpha|\gg 0.\]
\end{defn}
We denote by $ C_M(I^n)$ the set of all $C^\infty$ functions on $I^n$ which are in the class $M$.
\begin{rem}
Let $ M_1(t)= c r^t M(t)$, where $ c>0$, $r > 0$. We easily see that a function $ f\in\mathcal{E}(I^n)$ is in the class $M$ if and only if $f$ is in the class $M_1$ i.e. $ C_M(I^n) = C_{M_1}(I^n)$. We can then see that the class does not change when the function $ t\to m(t)$ is changed by the function
$ t\to m(t)+ at +b$ for some constants $a,b\in\R$. We will then suppose in the following that $ m(0)=0$.
\end{rem}
\begin{lem}
For all $ q\in\N$, there exist $C_q > 0$, $\rho_q > 0$ such that
\[ M(p+q)\leq C_q \rho_q ^p M(p), \,\,\,\,\forall p \in \N, \,\,\,\,p \gg 0.\]
\end{lem}
\begin{pre}
We have
\[ m(p+q) - m(p) = q m'(\theta), \,\,\,\,\mbox{where}\,\,\,\,\theta \in ]p, p+q[.\]
Since $m''(t) \leq \delta $, there exists $C >0$ such that $ m'(t) \leq \delta t + C$. We have then
\[ m(p+q) - m(p) \leq \delta p + (C + \delta q).\]
We put $ \rho_q = e^{q\delta}$ and $ C_q = e^{ q(q\delta  + C)}$, which proves the lemma.\hspace{5cm}$\square$
\end{pre}
\begin{lem}
$ C_M(I^n)$ is an algebra closed under differentiation.
\end{lem}
\begin{pre}
The function $ t\mapsto m(t)$ is convex and $ m(0)=0 $. For $ 0\leq q < p$, we have
\[ m(p-q)\leq \frac{p-q}{p} m(p)\,\,\,\,\mbox{and}\,\,\,\, m(q) \leq \frac{q}{p}m(p).\]
We obtain $ m(q) + m( p-q) \leq m( p)$, hence  $ M(q) . M( p-q) \leq M( p)$.\\
Using the last  inequality and the Leibniz formula, we deduce the first statement of the lemma.\\ The second statement follows immediately from lemma 3.3, or see \cite{Childress} (2.6).\hspace{3cm}$\square$
\end{pre}
\begin{rem}
If $ t \gg 0$, let  $ m(t) = t \log t$ i.e. $M(t) = t^t$, we have the analytic class. In the following we will consider the function $ t\mapsto m(t)$ such
the algebra $ C_M(I^n)$ contains strictly the analytic class. We therefore take the function $ t\mapsto m(t)$ of the form
\[ m(t)= t\log t + t\mu(t),\,\,\,\,\,\,\mbox{for}\,\,\,\,\,\,t\gg0\]
where $ \mu : [ 0, \infty[\to \R$ is an increasing function such that $\lim\limits_{t\to \infty}\mu(t) =\infty$.\\
In order to have $ m''(t) \leq \delta$ for $ t \gg 0$, we must suppose that $\mu (t) \leq a t $ for $ t \gg 0$, where $ a > 0$. We also suppose that the function
$ t\mapsto \mu (t)$ is in a Hardy field.
\end{rem}
\begin{prop}
$ C_M(I^n)$ is  closed under composition.
\end{prop}
\begin{pre}
We remark that the function $ t\mapsto t\mu (t)$ is convex, hence the proposition follows from \cite{Childress}, (2.5).
\end{pre}
\section{Quasianalyticity}
Let $ t\mapsto M(t)$ be as above, recall that
\[ M(t) = e^{m(t)} \,\,\,\,\,\, \mbox{and}\,\,\,\,\,\,m(t)= t \log t + t \mu (t)\,\,\,\,\,\,\mbox{for}\,\,\,\,\,\,t\gg0.\]
For each $ s\in\R^*_+$, we put
\[ \Lambda (s)= \inf\limits_{ t\geq t_0}M(t)s^{-t}\,\,\,\,\,\,\mbox{and}\,\,\,\,\,\, \lambda (s)= \inf\limits_{n\in\N,\,\, n\geq t_0}M(n)s^{-n},\]
where $t_0 \in \R_+$.\\
The minimum of the function $ t\mapsto M(t)s^{-t}$ is reached at a point $t$ where $ m'(t) = \log s$, and this point is unique, since the function $ t\mapsto m'(t)$ is increasing and $ \lim\limits_{t\to \infty}m'(t) = \infty$.\\
We define a function $ s\mapsto \omega(s)$ by
\[ \Lambda(s)= e^{-\omega(s)}.\]
We have the following relations
$$
\left\lbrace
\begin{array}{lll}
s  &  = &  e^{m'(t)}\\
\omega(s) & = & tm'(t) -m(t)
\end{array}
\right.
$$
or
$$
\left\lbrace
\begin{array}{lll}
s  &  = & e t  e^{\mu (t)+ t \mu' (t)}\\
\omega(s) & = & t + t^2 \mu' (t)
\end{array}
\right.
$$
We can easily invert the previous system, so we  obtain the following system:
$$
\left\lbrace
\begin{array}{lll}
t  &  = & s\omega'(s)\\
m(t)  & = & s\omega'(s)\log s - \omega(s)
\end{array}
\right.
$$
Using the fact that $ m(t)= t\log t + t\mu(t)$, we obtain
$$
\left\lbrace
\begin{array}{lll}
t  &  = & s\omega'(s)\\
-\mu(t)  & = & \log \omega'(s) + \frac{ \omega(s)}{ s \omega'(s)}.
\end{array}
\right.
$$
We can then see that the function $ s \mapsto \omega(s)$ is increasing and when $ t \to \infty$, we have
\[ s\omega'(s) \to \infty \,\,\,\,\,\, \mbox{and}\,\,\,\,\,\,\log \omega'(s) + \frac{ \omega(s)}{ s \omega'(s)} \to -\infty.\]
Thus $\omega'(s)\to 0$ when $ s\to \infty$.
\begin{lem}
For $ s \gg 0$, we have
\[ e^{-\delta} \lambda(s) \leq \Lambda(s) \leq \lambda(s).\]
\end{lem}
\begin{pre}
Put $\alpha(t) = m(t) - t \log s$, if $ \alpha'(t_0) = 0$, we have $\Lambda(s) = e^{\alpha(t_0)}$. Let $ n_0\in \N$ such that $ |n_0 - t_0| < 1$, hence
$\lambda(s) = e^{\alpha(n_0)}$.\\
Note that
\[ \alpha(n_0) - \alpha(t_0) = (n_0 - t_0)\alpha'((1-\theta )n_0 + \theta t_0 ),\,\,\,\,\,\mbox{where}\,\,\,\,\, 0 < \theta <1.\]
Since $ \alpha''(t) \leq \delta$ and $\alpha't_0)= 0$, we get  $e^{-\delta} \lambda(s)$. The second inequality is trivial.
\hspace{12cm}$\square$
\end{pre}
\begin{prop}
The following conditions are equivalent:
\begin{enumerate}
\item [ a)] $ \sum\limits_{p =0}^\infty \frac{ M(p)}{ M(p+1)} = \infty$,
\item [ b)] $\int_{s_0}^\infty\frac{ \omega(s)}{s^2} ds = \infty$,  for some $s_0$,
 \item [ c)]$\int_{s_0}^\infty\frac{ \log \lambda(s)}{s^2} ds = - \infty$,  for some $s_0$.
\end{enumerate}
\end{prop}
\begin{pre}
If $ p\in \N$, we have
\[ m'(p) \leq m(p+1) - m(p)\leq m'(p+1),\]
hence
\[ \sum\limits_{p =0}^\infty \frac{ M(p)}{ M(p+1)} = \infty \Leftrightarrow \int_{t_0}^\infty e^{-m'(t)}dt = \infty .\]
From above
\[ \int_{t_0}^\infty e^{-m'(t)}dt = \int_{s_0}^\infty \frac{d( s\omega'(s))}{s} ds,\]
since $ \omega'(s) \to 0$ when $ s \to \infty$, we get
\[ \int_{s_0}^\infty \frac{d( s\omega'(s))}{s} ds = \infty \Leftrightarrow \int_{s_0}^\infty\frac{ \omega(s)}{s^2} ds = \infty,\]
which proves $ a) \Leftrightarrow b)$. By lemma 4.1, we deduce the equivalence
 $ b) \Leftrightarrow c)$.
\end{pre}
\begin{defn}
We say that the algebra $ C_M(I^n)$ is quasianalytic if any $ f\in C_M(I^n)$ and any $ x\in I^n$, the Taylor series $T_xf$ of $f$ at $x$ determines uniquely $f$ near $x$.
\end{defn}
\noindent      If $C_M(I^n)$ is quasianalytic, we say that the class $M$ is quasianalytic.\\
A well-known result of Denjoy-Carleman gives a characterization of quasianalyticity in terms of the function $ t\mapsto M (t)$.
\begin{equation}
C_M(I^n)\,\,\,\,\,\mbox{ is quasianalytic } \,\,\,\,\,         \Leftrightarrow           \,\,\,\,\,\sum\limits_{p =0}^\infty \frac{ M(p)}{ M(p+1)} = \infty.
\end{equation}
\begin{rem}
As  a consequence of Proposition 4.3, we see that if the class $M$ is quasianalytic, then $\lim\limits_{s\to \infty}
\frac{\omega(s)}{s^q}= \infty$,  for each $0\leq q < 1$.
\end{rem}
{\bf{Probably the converse of this statement is true.}}\\
Let us examine this fact if the class M is analytic.\\
If the class $M$ is analytic i.e. $m(t)= t\log t$,  we have  $ \omega(s) = s \omega'(s)$, hence $ \omega(s) = C s$, for some $ C\in\R$.\\
The converse is also true.
\begin{thm}
If $ \omega(s) \sim s$ when $ s\to \infty$, then any $ f\in C_M(I^n)$ is analytic.
\end{thm}
\begin{pre}
We have $\lim\limits_{s\to \infty}\frac{\omega(s)}{s} =1$, hence there exist $ a > 0$, $ c > 0$, such that
\[ \omega(s) \geq c s,\,\,\,\,\,\,\forall s \geq a,\]
thus, we have
\[e^{- \omega(s)}  \leq  \frac{q !}{c^q s^q},\,\,\,\,\,\,\forall q \in \N,\,\,\,\,\,\,\forall s \geq a .\eqno{(*)}\]
Since $ m'(t)\to  \infty $  when $ t\to \infty$, there exists $ N_0\in \N$ such that $$ e^{m'(t)} \geq c,\,\,\,\,\,\, \,\,\,\,\,\, \forall  t > N_0.$$
Recall that $$\lambda (s)= \inf\limits_{n\in\N,\,\, n\geq t_0}M(n)s^{-n},$$ for some $ t_0\in \R_+$. Note that we can suppose $ N_0 > t_0 + 1$.\\
We will now show that
\[ M(p) \leq e^\delta c^{-p}p!,\,\,\,\,\,\, \,\,\,\,\,\, \forall  p > N_0,\,\,\,\, p\in \N .\]
For $  p > N_0,\,\,\,\, p\in \N$, let $ r\in \R$ be such that $ |p-r| < \frac{1}{2}$, and put $ s = e^{m'(r)}$.
We have then
\[ s \geq a,\,\,\,\,\,\, \,\,\,\,\,\, \Lambda(s)= M(r) s^{-r},\,\,\,\,\,\, \,\,\,\,\,\,\lambda(s)= M(p) s^{-p}.\]
By lemma 4.1 and $(*)$, we have
\[ \lambda(s)= M(p) s^{-p} \leq e^\delta \Lambda(s)\leq e^\delta c^{-p} p! s^{-p},\]
hence, for each $p\in\N$, $ p > N_0$, we have
\[ M(p) \leq e^\delta c^{-p}p!,\]
which proves the theorem.\hspace{11cm}$\square$
\end{pre}
\begin{exa}
Let $\mu(t)= \log\log t$, i.e. $m(t) = t \log t + t \log\log t$.
The class $M(t) = e^{m(t)}$ is quasianalytic. Indeed,  we have
\[ s = e^{m'(t)}  = e t e^{\frac{1}{\log t}}\log t \sim e t \log t .\]
On the other hand
\[ \omega(s) = t m'(t) - m(t) = t + \frac{t}{\log t} \sim t.\]
Since $ \log s = m'(t)$ and $ m'(t) \sim \log t$, we have
\[ \frac{\omega(s)}{s^2} \sim  \frac{1}{e s \log s},\,\,\,\,\,\,\,\mbox{when} \,\,\,\,\,\,\,s\to\infty.\]
By $b)$ Proposition 4.3., we deduce that the class $M$ is quasianalytic.
\end{exa}
\section{Shifted Denjoy-Carleman classes}
From a  class $M$, one can construct other classes as follows: \\ Let $ M(t) = e^{m(t)}$.
For each $ p\in\N^*$, we put $$ m_p(t) = m(pt) \,\,\,\,\,\,\,\mbox{and} \,\,\,\,\,\,\,M_p(t) = e^{m_p(t)}.$$ 
It is clear that $ C_M(I^n) \subset C_{M_p}(I^n),\,\,\,\,\forall p\in \N$.
It is easy to see that the function $$m_p : t\mapsto m_p(t) = m(pt)$$ satisfies
the properties 1), 2), and 3) of section 2. We deduce that  Lemma 3.5. and Proposition 3.8. are still true for the rings $ C_{M_p}(I^n)$.
\begin{rem}
However, it can happen that a class $M $ is quasianalytic but the class $M_p$    is not quasianalytic for each  $p\in \N,\,\,\,\, p >1$.
\end{rem}
\begin{exa}
We have already seen that if $m(t) = t \log t + t \log\log t$, then the class $ M(t)= e^{m(t)}$ is quasianalytic. \\Let us see that the class $ M_p$ is
not quasianalytic for each $p\in \N,\,\,\,\, p >1$. We have
\[ m_p(t) = m(pt)= p t \log pt + pt \log\log p t.\]
By Remark 3.2, we can take $ m_p(t) = p t \log t + pt \log\log pt$  without changing the class.\\
We put
\[ \Lambda_p(s) = \inf\limits_{t\geq t_0} M_p(t)s^{-t},\,\,\,\,\,\,\mbox{for some}\,\,\,\, t_0 >0\,\,\,\,\,\,\mbox{and}\,\,\,\,\,\,
e^{- \omega_p(s)} = \Lambda_p(s).\]
By proceeding as in Example 4.9, we find
\[ \frac{\omega_p(s)}{s^2} \sim \frac{ p^2}{e} \frac{1}{ s^{\frac{2p-1}{p}}\log s}.\]
We see that if $ p> 1$, then
 $ \int_{s_0}^\infty \frac{\omega_p(s)}{s^2} < \infty$. By Proposition 4.3, the class $M_p$ is not quasianalytic.\\
 In the above proof It may be sufficient to show that $M_2$ is not quasianalytic since \\$ C_{M_2}(I^n)\subset C_{M_p}(I^n),\,\,\,\forall p \geq 2$.
\end{exa}
Here is a class $ \tilde{M}$ where this fact can not happen. 
This class   will allow us to build a non-noetherian quasianalytic system.
Let
\[ \tilde{m}(t) = t \log t + t \log \log \log t,\,\,\,\,\,\mbox{and}\,\,\,\,\,\forall p\in \N^*, \,\,\tilde{m}_p(t)= \tilde{m}(pt),\,\,\,  t > e^e. \]
We claim that for each $p\in \N^*$, the class $\tilde{M}_p(t)= e^{\tilde{m}_p(t)} $ is quasianalytic. Indeed, it is easy to see that

\[ C_{\tilde{M}_p}(I^n) \subset C_{M}(I^n),\,\,\,\forall p\in \N^*,\]
where $ M(t) = e^{m(t)}$, with $m(t) = t \log t + t \log\log t$. We deduce, by (4.1),  that the class $\tilde{M}_p$ is quasianalytic for each $ p\in \N^*$.\\
We denote by $ C_{\tilde{M}^\infty}(I^n)$ the ring $\bigcup\limits_{p=1}^\infty C_{\tilde{M}_p}(I^n)$.
Note that the classes $\tilde{M}_p$ strictly contain the analytic class.\\
We denote by $ C_{\tilde{M}^\infty}(n)$ the ring of germs at the origin in $\R^n$ of functions from $ C_{\tilde{M}^\infty}(I_\epsilon^n)$, for some small $\epsilon >0$.
It is clear that $ C_{\tilde{M}^\infty}(n)$ is a local ring with maximal ideal $ \underline{m}_{\tilde{M}^\infty}(n)$
generated by $\{ x_1,\ldots, x_n\}$ and its completion with the $ \underline{m}_{\tilde{M}^\infty}$-topology is the ring of formal power series
$\R[[X_1,\ldots, X_n]]$.
We consider the system
\[ \mathcal{C}_{\tilde{M}^\infty} =\{ C_{\tilde{M}^\infty}(n),\,\,n\in \N^*\}.\]
It can easily be seen that the system $\mathcal{C}_{\tilde{M}^\infty}$ contain strictly the system of the analytic germs and is closed under
composition, partial differentiation and implicit function. We also note that for each $n\in\N$, the Taylor expansion at the origin map,
$ f\mapsto \hat{f}$, is injective on  $C_{\tilde{M}^\infty}(n)$.
\section{ Well behaved quasianalytic system}
In this section we fix the quasianalytic system
\[ \mathcal{C}_{\tilde{M}^\infty} =\{ C_{\tilde{M}^\infty}(n),\,\,n\in \N^*\},\] constructed from the classes
\[ \tilde{M}_p(t) = e^{\tilde{m}_p(t)},\,\,\,\,\,\,\mbox{where}\,\,\,\,\,\, \tilde{m}_p(t) = \tilde{m}(pt)= p t \log pt + pt \log\log p t.\]
We are concerned here with local homomorphisms
\[ \Phi: C_{\tilde{M}^\infty}(n)\to C_{\tilde{M}^\infty}(k),\]
i.e. homomorphisms such that $\Phi( \underline{m}_{\tilde{M}^\infty}(n))\subset \underline{m}_{\tilde{M}^\infty}(k)$. Since $\Phi$ is a local homomorphism,
we consider its natural extension to the completion $\hat{\Phi}: \R[[ X_1,\ldots, X_n]] \to  \R[[ X_1,\ldots, X_k]]$.\\
Let \[ \hat{\Phi}_* : \frac{\R[[ X_1,\ldots, X_n]]}{C_{\tilde{M}^\infty}(n)}\to \frac{\R[[ X_1,\ldots, X_k]]}{C_{\tilde{M}^\infty}(k)}\]
be the homomorphism of groups induced by $\hat{\Phi}$ in the obvious manner.
\begin{defn}
We say that the homomorphism $\Phi$ is {\it{ strongly injective}}, if the homomorphism $\hat{\Phi}_*$ is injective.
\end{defn}
$\Phi$ is strongly injective  can be expressed as follows:\\
 \[\mbox{if}\,\,\,\,\, \hat{f}\in  \R[[ X_1,\ldots, X_n]]\,\,\,\,\,\mbox{ is such that}\,\,\,\,\,  \hat{\Phi}(\hat{f}) \in C_{\tilde{M}^\infty}(k),
 \,\,\,\,\,\mbox{ then}\,\,\,\,\, \hat{f}\in C_{\tilde{M}^\infty}(n).\]
We consider the local homomorphisms
\[ e: C_{\tilde{M}^\infty}(n) \to C_{\tilde{M}^\infty}(n),\,\,\,\,\,\,\,r_d: C_{\tilde{M}^\infty}(n) \to C_{\tilde{M}^\infty}(n)\,\,\,\,\,\,\,\mbox{where}\,\,\,\,\,\,
d\in\N,\] defined by:
\[ e(f)(x_1,x_2,\ldots,x_n) = f(x_1x_2,x_2,\ldots,x_n),\,\,\,\,\,\,\, r_d(f)(x_1,x_2,\ldots,x_n)= f(x_1^d,x_2,\ldots, x_n).\]
\begin{defn}
A quasianalytic system $\mathcal{C}=\{ \mathcal{C}_n\,/\,n\in\N\}$ is called {\it{ well behaved}}, if $e$ and for each $d\in\N$, $r_d$ are strongly injective.
\end{defn}
\begin{exa}
If for each $n\in\N$, $\mathcal{C}_n$ is the ring of germs, at the origin in $\R^n$, of real analytic functions, the obtained system is well behaved.\\ 
If for each $n\in\N$, $\mathcal{C}_n$ is the ring of germs, at the origin in $\R^n$, of Nash functions, the obtained system is well behaved, see \cite{Elkhadiri}
\end{exa}
\begin{prop}
The quasianalytic system $\mathcal{C}_{\tilde{M}^\infty} =\{ C_{\tilde{M}^\infty}(n),\,\,n\in \N^*\}$ is well behaved.
\end{prop}
\begin{pre}
We must show  that the local homomorphisms  $e$ and $ r_d$ are strongly injective.
We will give the proof for $ r_d$, the proof for $e$ is the same.\\
Let $\hat{f}\in  \R[[ X_1,\ldots, X_n]]$ be such that $ \hat{r_d}(\hat{f})= \hat{g} \in C_{\tilde{M}^\infty}(n)$. There exist $p\in\N$, $ \epsilon > 0$
 such that $ g\in C_{\tilde{M}_p}(I_\epsilon^n)$. By \cite{Bierston}, Theorem 4.1, there is a mapping $ I_\epsilon^n\to \R[[ X_1,\ldots, X_n]]$,
$ a\mapsto H^*_{r_d(a)}$, such that:
\[ \forall a\in I_\epsilon^n,\,\,\,\,\, \hat{g}_a = H^*_{r_d(a)}\circ \hat{r_d}_a.\]
Where $ \hat{g}_a$ is the Taylor expansion of the function $g$ at the point $a$ and if \\
$\sigma = (\sigma_1,\ldots,\sigma_n): U\subset \R^n\to \R^n$ is a $C^\infty$ mapping and
$$ H^*_{r_d(a)} = \sum\limits_{\omega=(\omega_1,\ldots,\omega_n)\in \N^n} b_\omega(a) X_1^{\omega_1}\ldots X_n^{\omega_n}$$
a formal series. The formal series $H^*_{\sigma(a)}\circ \hat{\sigma}_a $ is given by:
$$ H^*_{\sigma(a)}\circ \hat{\sigma}_a  = \sum\limits_{\omega=(\omega_1,\ldots,\omega_n)\in \N^n} b_\omega(a) [ \hat{( \sigma_1 - \sigma_1(a))}]^{\omega_1}
\ldots [ \hat{( \sigma_n - \sigma_n(a))}]^{\omega_n},$$
where  $\hat{( \sigma_j - \sigma_j(a))}, j= 1\ldots,n$, is the Taylor expansion of the function $ x\mapsto \sigma_j(x) - \sigma(a)$ at the point $a$.\\
Let $ U$ be an open neighborhood of $I_\epsilon^n$ in $\R^n$,
the mapping $ r_d: U\subset \R^n \to \R^n$ is analytic, proper and generically a submersion. By  Glaeser's composite function  Theorem, \cite{Glaeser} Theorem II,
there exists a $C^\infty$ function, $h$, on $U$ such that $ g = h\circ r_d$ on $U$. Since $ r_d(I_\epsilon^n) = I_\epsilon^n$, by \cite{Bierston}, lemma 3.1, we see that the restriction of the function $h$ to $I_\epsilon^n$ is in
 $  C_{\tilde{M}_{dp}}(I_\epsilon^n)\subset C_{\tilde{M}^\infty}(I_\epsilon^n)$, which proves the proposition.\\
 For the mapping $e$ we use the fact that $\overline{ e(I_\epsilon^n)} = I_\epsilon^n$, and instead of using \cite{Bierston}, lemma 3.1, we use \cite{Bierston}, lemma 3.4.
\end{pre}
Let us recall the following result proved  in \cite{Elkhadiri}, Corollary 5.8.
\begin{prop}
Every well behaved noetherian system is contained in the analytic system.
\end{prop}
We remark that every noetherian system is a quasianalytic system.
\begin{cor}
The system \[ \mathcal{C}_{\tilde{M}^\infty} =\{ C_{\tilde{M}^\infty}(n),\,\,n\in \N^*\}\] is not noetherian, i.e. there exists $ m\in\N$, $ m > 1$, such that the ring $ C_{\tilde{M}^\infty}(m)$ is not noetherian.
\end{cor}
\begin{pre}
Since the system $\mathcal{C}_{\tilde{M}^\infty}$ contains strictly the analytic system, the corollary is 
 is a direct consequence of Proposition 6.4. and Proposition 6.6.
\end{pre}
\section{Other examples of non-Noetherian quasianalytic rings}
Fix a polynomially bounded o-minimal structure $\mathcal{R}$ on the field of reals and denote by $\mathcal{D}_n$ the ring of those
quasianalytic function germs at $0\in \R^n$ which are definable in $\mathcal{R}$. Suppose that  the system $\mathcal{D} = \{ \mathcal{D}_n\,/\, n\in\N^*\}$ contains
 strictly the  analytic system.
\begin{thm}
The system $\mathcal{D} = \{ \mathcal{D}_n\,/\, n\in\N^*\}$ is not noetherian, i.e. there exists $ m\in\N,\,\,\, m >1$, such that the ring $\mathcal{D}_m$ is not noetherian.
\end{thm}
\begin{pre}
From the above, it remains to show that the mappings $r_d: \mathcal{D}_n\to \mathcal{D}_n$ and $e: \mathcal{D}_n\to \mathcal{D}_n$   are strongly injective.\\
By following the proof of the Proposition 6.4. we note that the function $h$ given by  Glaeser's composite function  Theorem is definable on $I_\epsilon^n$, since $r_d(I_\epsilon^n) = I_\epsilon^n$. Here  one does not need to evoke the result of \cite{Bierston}, lemma 3.1. The same thing is true for the mapping $e$, since $\overline{ e(I_\epsilon^n)} = I_\epsilon^n$.
\end{pre}


\begin{thebibliography}{5}
\bibitem{Broglia} \bf F. Acquistapace, F. Broglia, M. Bronshtein, A. Nicoara, and N. Zobin, {\em Failure of the Weierstrass preparation theorem in quasianalytic local rings }, Eprint arXiv: 1212.4265, 2012.
\bibitem{Bierston} \bf A. Belotto da Silva, I. Biborski, and E. Bierstone, {\em Solutions of quasianalytic equations}, arXiv: 1605.01425 [math].
\bibitem{Carleman} \bf T. Carleman, {\em Les  fonctions quasi-analytiques }, Gauthiers Villars Paris  (1926)
\bibitem{Childress} \bf C.L.Childress, {\em Weierstrass division in quasianalytic local ring}, Can. J. Math., Vol.XXVIII.
\bibitem{Denjoy} \bf A. Denjoy, {\em Sur les fonctions quasianalytiques de la variable r\'eelle}, C. R. Acad. Sci. Paris (123) (1921).
\bibitem{Elkhadiri} \bf A. Elkhadiri, {\em Link between Noetherianity and the Weierstrass division theoremon some quasianalytic local rings}, 
Proc. Amer. Math. Soc. 140 (2012), 3883-3892.
\bibitem{Glaeser} \bf G. Glaeser, {\em  fonctions composées différentiables}, Ann. of Math. Vol. 77. N0. 1 Januar, 1963.
\bibitem{Komatsu} \bf H. Komatsu, {\em The implicit function theorem for ultradifferentiable mapping }, Proc. Japan Acad. Ser. A Math. Sci. 55 (1979), 69-72.
\bibitem{Miller} \bf C.Miller, {\em Infinite diffentiability in polynomially bounded o-minimal structure}, Proc. Amer. Math. Soc. 123 (1995), 2551-2555.
\bibitem{Parusinski} \bf A. Parusinski and J. P. Rolin, {\em A note on the Weierstrass preparation theorem in quasianalytic local rings}, Canad. Math. Bull. Vol. 57 (4) 2014 pp. 614-620.






\end{thebibliography}
\end{document}